\documentclass[10pt, reqno]{amsart}
\UseRawInputEncoding
\usepackage{amsmath,amsfonts,amssymb,color,a4,epsfig,graphics}
\usepackage[colorlinks=true,hyperindex=true]{hyperref}
\usepackage[english]{babel}
\usepackage{enumerate}
\usepackage{xcolor}
\allowdisplaybreaks
\RequirePackage{ifthen}
\usepackage{tikz}
\def\build#1_#2^#3{\mathrel{\mathop{\kern 0pt#1}\limits_{#2}^{#3}}}
\usetikzlibrary{arrows}

\newcommand{\R}{{\mathbb{R}}}

\newcommand{\Z}{{\mathbb{Z}}}
\newcommand{\N}{{\mathbb{N}}}

\newcommand{\Bc}{\mathcal{B}}
\newcommand{\Cc}{\mathcal{C}}
\newcommand{\Dc}{\mathcal{D}}
\newcommand{\Ec}{\mathcal{E}}
\newcommand{\Fc}{\mathcal{F}}
\newcommand{\Gc}{\mathcal{G}}
\newcommand{\Hc}{\mathcal{H}}
\newcommand{\Ic}{\mathcal{I}}

\newcommand{\Nc}{\mathcal{N}}
\newcommand{\Sc}{\mathcal{S}}
\newcommand{\Vc}{\mathcal{V}}

\def\rmi{{\rm i}}
\def\rmadj{{\rm adj}}

\newtheorem{theorem}{Theorem}[section]
\newtheorem{proposition}[theorem]{Proposition}
\newtheorem{lemma}[theorem]{Lemma}
\newtheorem{remark}[theorem]{Remark}
\newtheorem{definition}[theorem]{Definition}
\newtheorem{corollary}[theorem]{Corollary}

\begin{document}
\title[\tiny Limiting absorption principle for long-range perturbation in the discrete...]{Limiting absorption principle for long-range perturbation in the discrete triangular lattice setting}
\author[\tiny Nassim Athmouni]{Nassim Athmouni$^{1}$}
\address{$^{1}$Universit\'e de Gafsa, Campus Universitaire, 2112, Tunisie}
\email{\tt nassim.athmouni@fsgf.u-gafsa.tn}
 \author[\tiny Marwa Ennaceur]{Marwa Ennaceur$^{2}$}
 \address{$^{2}$Department of Mathematics, College of Science, University of Hail, Ha’il 81451, \hspace*{0,5cm}Saudi Arabia}
\email{\tt mar.ennaceur@uoh.edu.sa}
\author[\tiny Sylvain Gol\'enia]{Sylvain Gol\'enia$^{3}$}
\address{$^{3}$Univ. Bordeaux,
Bordeaux INP, CNRS, IMB, UMR 5251, F-33400 Talence, France}
\email{\tt sylvain.golenia@math.u-bordeaux.fr}
\author[\tiny Amel Jadlaoui]{Amel Jadlaoui$^{4}$ }
 \address{$^{4}$Universit\'e de Sfax. Route de la Soukra km 4 - B.P. n° 802 - 3038 Sfax, Tunisie}
\email{\tt amel.jadlaoui@yahoo.com}

\subjclass[2010]{81Q10, 47B25, 47A10, 05C63,  08B15}
\keywords {Commutator, Mourre estimate, Limiting Absorption Principle, Discrete Laplacian, Triangular lattice}
\begin{abstract} We examine the discrete Laplacian acting on a triangular lattice, introducing long-range perturbations to both the metric and the potential. Our goal is to establish a Limiting Absorption Principle away from
possible embedded eigenvalues. Our study relies on a positive commutator technique.
\end{abstract}
\maketitle
\tableofcontents
\section{Introduction and main result}
In recent years, spectral graph theory has attracted significant attention, particularly in the study of various for different types of discrete Laplacians \cite{AD,AEG,BGJ,S,GG,Ch,Mic,Gk,AT,BBJ}
and their magnetic analogs \cite{GT,BGKLM,GoMo,ABDE,PR1,HiSh}. One approach to analyzing the essential spectrum of these operators is based on 
a positive commutator technique. For instance, the authors in \cite{S,BoSa} study the case of $\Z^d$, while \cite{AF} and \cite{GG} analyze binary trees. similarly, 
\cite{MRT} investigates a general family of graphs and \cite{AEG} works on a discrete version of cusps and funnels. 

In \cite{PR},
the authors study the spectral theory of Schr{\"o}dinger
operators acting on perturbed periodic discrete graphs. They consider two types of perturbations: a long-range potential and a short-range modification of the metric. 
Using the Mourre estimate and take advantage of a Floquet-Bloch decomposition, they prove a Limiting Absorption Principle. 
In the present work, we focus on a specific case: the triangular lattice (see Figure 1). Our goal is to obtain similar spectral results but for a broader class of perturbations. We introduce long-range perturbation to both the potential and the metric. Our approach relies on a Mourre estimate technique.

We begin with some standard definitions from graph theory. An infinite, connected \textit{graph} $\Gc$ is a triplet $(\Vc,\Ec,m)$, where $\Vc$ is the countable set of \textit{vertices}, $m : \Vc \rightarrow (0, \infty)$ is a \textit{weight} and $\Ec:\Vc\times\Vc \rightarrow  [0,+\infty)$ (the edges) is symmetric.
Given two vertices $n$ and $l$, we say that $n$ and $l$ are \textit{neighbors} if $\Ec(n,l)>0$. We denote this relationship by $n\sim l$. The set of neighbors of $n$ is denoted
by $\Nc_n$. 
We denote by $\Cc(\Vc):=\{f:\Vc\longrightarrow\mathbb{C}\}$ the space of complex-valued functions acting on the set of vertices $\Vc$.
Now, we consider the Hilbert space: \[\ell^2(\Vc,m):=\left\{ f\in \Cc(\Vc); \sum_{n\in\Vc} m(n)|f(n)|^2 <\infty \right\},\]
equipped with the scalar product, $\langle f,g\rangle:=\displaystyle \sum_{n\in \Vc} m(n)\overline{f(n)}g(n)$.

Now, we define our model as follows. Set \[\Vc:=\Big\{\sum_{j=1}^2k_j v_j\hbox{; } k:=(k_1,k_2)\in \mathbb{Z}^2\Big\}\hbox{, } v_1:=(1,0)\hbox{, }v_2:=(\frac{1}{2},\frac{\sqrt{3}}{2}).\]
We define \begin{align*}
\Ec: \Vc\times\Vc&\to \{0,1\}
\\
(l,n)&\mapsto\Ec(l,n):=\left\{
                        \begin{array}{ll}
                          1, & \hbox{if } |l-n|_{\mathbb{R}^2}=1; \\
                          0,
                        \end{array}
                      \right.
\end{align*}
where $|n|_{\mathbb{R}^2}:=\sqrt{n_1^2+n_2^2}$. We introduce
\begin{align}\label{Na}\Nc_{n}:=\Big\{l\in\Vc\hbox{; }|l-n|_{\mathbb{R}^2}=1,\ n\in\Vc \Big\}=\Big\{n\pm v_{1}\hbox{, }n\pm v_{2}\hbox{, }n\pm v_3 \Big\},\end{align}
where $ v_3:=v_{1}-v_2$.
We note that $\sharp \Nc_n=6$, for all $n\in\Vc.$

\begin{figure}[htbp]
  \centering
  \begin{tikzpicture}[x=0.5cm, y=0.5cm,line width=1pt, scale=1,transform shape]
   \draw (0, 3) -- (2, 0);
   \draw (0, 9) -- (6, 0);
   \draw (0, 15) -- (10, 0);
   \draw (4, 15) -- (14, 0);
   \draw (8, 15) -- (18, 0);
   \draw (12, 15) -- (20, 3);
   \draw (16, 15) -- (20, 9);
   \draw (4, 15) -- (0, 9);
   \draw (8, 15) -- (0, 3);
   \draw (12, 15) -- (2, 0);
   \draw (16, 15) -- (6, 0);
   \draw (20, 15) -- (10, 0);
   \draw (20, 9) -- (14, 0);
   \draw (20, 3) -- (18, 0);
   
   \draw (-2, 12) -- (22, 12);
   \draw (0, 9) -- (20, 9);
   \draw (-2, 6) -- (22, 6);
   \draw (0, 3) -- (20, 3);
   \draw (0, 15) -- (20, 15);
   \draw (2, 0) -- (18, 0);
    \draw[dashed] (20, 9) -- (22, 9);
    \draw[dashed] (0, 9) -- (-2, 9);
    \draw[dashed] (20, 3) -- (22, 3);
    \draw[dashed] (0, 3) -- (-2, 3);
   \draw[dashed] (6, 0) -- (6, -2);
   \draw[dashed] (14, 0) -- (14,-2);
   \draw[dashed] (4, 15) -- (4, 17);
   \draw[dashed] (16, 15) -- (16,17);
   \foreach \i/\j in {2/-2,6/-1,10/0,14/1,18/2}{
	\node[draw=none, circle, fill=blue, minimum size=1mm, label={[font=\small]above right:(\j, $\sqrt{3}$)}] at (\i,12){} ;}

  \foreach \i/\j in {4/-3,8/-1,12/1,16/3}{
	\node[draw=none, circle, fill=blue, minimum size=1mm, label={[font=\small]above right:($\frac{\j}{2},\frac{\sqrt{3}}{2}$)}] at (\i,9){} ;}

  \foreach \i/\j in {2/-2,6/-1,10/0,14/1,18/2}{
	\node[draw=none, circle, fill=blue, minimum size=1mm, label={[font=\small]above right:(\j,0)}] at (\i,6){} ;}

  \foreach \i/\j in {4/-3,8/-1,12/1,16/3}{
	\node[draw=none, circle, fill=blue, minimum size=1mm, label={[font=\small]above right:($\frac{\j}{2},-\frac{\sqrt{3}}{2}$)}] at (\i,3){} ;}

  \end{tikzpicture}
  \caption{Triangular lattice}
\end{figure}
In this sequel, we often identify the vertices of $\Vc$ with $\mathbb{Z}^2$, using the canonical map
\begin{align}\label{map}
\nonumber \mathbb{Z}^2 &\to \Vc
\\
(k_1,k_2)&\mapsto k_1v_1+k_2v_2.
\end{align}
We define the Laplacian $\Delta_m$ by
\begin{align}\label{dm}
\nonumber\Delta_{m}: \ell^2(\Vc,m)&\to \ell^2(\Vc,m)
\\
f&\mapsto\Delta_{m}f(n):=\frac{1}{6m(n)}\sum_{i=1}^3f(n+v_i)+f(n-v_i).
\end{align}
Let $\eta$ be a real-valued function on $\Vc$ such that:  \begin{align*}
 (H_0)\quad& m(n):=(1+\eta(n)),\
\mbox {and } \inf_n\eta(n)>-1,\ \eta(n)\rightarrow0  \mbox{ if } |n|\rightarrow\infty.\end{align*} 
$\Delta_m$ is bounded and self-adjoint on $\ell^2(\Vc;m)$. In the case where $\eta\equiv0$, one recovers the discrete Laplacian on the triangular lattice $\Delta_T$. It is known that its spectrum is $\left[-\frac{1}{2},1\right]$ and absolutely continuous (cf. \cite{AIM} and Lemma \ref{specter} below). Now, we seek similar results for $\Delta_m+V$, for long-range potentials $V.$ 

For a function $G: \Vc\rightarrow \mathbb{C}$, we denote by $G(Q_1,Q_2)$ the operator of multiplication by $G$. In particular, $(G(Q_1,Q_2)f)(n_1,n_2):=G(n_1,n_2)f(n_1,n_2)$, for all $f\in \Dc(G(Q_1,Q_2))$, where $$\Dc(G(Q_1,Q_2)):=\left\{f\in \ell^2(\Vc,m);\ n\mapsto G(n_1,n_2)f(n)\in \ell^2(\Vc,m)\right\}.$$
Let $V$ be a real-valued bounded function on $\Vc$, and $H_m:=\Delta_m+V(Q)$, such that: \begin{align*}(H'_0)\quad&  V(n)\to0 \mbox{ if } |n|\to\infty.\end{align*} 
Since $V(Q)$ is a compact operator, as uniform limit of finite rank operators given by $1_{\|\cdot\|_{\mathbb{R}^2}\leq R}V$, with $R\in\mathbb{N}$. The operator $H_m$ is bounded and self-adjoint on $\ell^2(\Vc,m)$. In fact, it is a kind of compact perturbation of $\Delta_{T}$, see Proposition \ref{compact}. Moreover, we have $\sigma_{\rm ess}(H_m)=\sigma_{\rm ess}(\Delta_T)$, where $\sigma_{\rm ess}(\cdot)$ denotes the essential spectrum, see Proposition \ref{com} for a precise statement.

Now, we aim for a more refined spectral property and ask for further decay. Let \begin{align}\label{lambda}\Lambda\left(n_1,n_2\right):=\langle n_1\rangle+\langle n_2\rangle,\end{align} where $\langle \cdot\rangle:=\sqrt{\frac{1}{2}+|\cdot|^2}$. Note that $\Lambda(Q_1,Q_2)$ is an unbounded self-adjoint operator. From now on,
\textbf{we fix} $\boldsymbol{\varepsilon>0}$ and introduce different hypotheses of decay for the metric: 
 \begin{align*}
 (H_1)\quad& \sup_{(n_1,n_2)\in\mathbb{Z}^2}\Lambda^{\varepsilon}(n_1,n_2)\langle n_1\rangle\left|\eta(n_1,n_2)-\eta(n_1+1,n_2)\right|<\infty,
 \\
 (H_2)\quad& \sup_{(n_1,n_2)\in\mathbb{Z}^2}\Lambda^{\varepsilon}(n_1,n_2)\langle n_2\rangle\left|\eta(n_1,n_2)-\eta(n_1,n_2+1)\right|<\infty,
 \\
 (H_3)\quad& \sup_{(n_1,n_2)\in\mathbb{Z}^2}\Lambda^{\varepsilon}(n_1,n_2)\langle n_1-n_2\rangle\left|\eta(n_1,n_2)-\eta(n_1+1,n_2-1)\right|<\infty.
 \end{align*}

Similarly, for the potential:
\begin{align*}
(H'_1)\quad&\sup_{(n_1,n_2)\in\mathbb{Z}^2}\Lambda^{\varepsilon}(n_1,n_2)\langle n_1\rangle\left|V(n_1,n_2)-V(n_1+1,n_2)\right|<\infty,\\
(H'_2)\quad&\sup_{(n_1,n_2)\in\mathbb{Z}^2}\Lambda^{\varepsilon}(n_1,n_2)\langle n_2\rangle\left|V(n_1,n_2)-V(n_1,n_2+1)\right|<\infty,\\
(H'_3)\quad&\sup_{(n_1,n_2)\in\mathbb{Z}^2}\Lambda^{\varepsilon}(n_1,n_2)\langle n_1-n_2\rangle\left|V(n_1,n_2)-V(n_1+1,n_2-1)\right|<\infty.
\end{align*}
Here, we have used the identification given by \eqref{map}.

Set $\kappa(H_m):=\left\{-\frac{1}{2},\ -\frac{1}{3},\ 1\right\}$. We denote by $\sigma_{\rm p}(\cdot)$ the set of pure point spectra.
We state our main theorem: 
\begin{theorem}\label{t:LAP} 
Suppose that $(H_0),\ (H_1),\ (H_2),\ (H_3),\ (H'_0),\ (H'_1),\ (H'_2)$ and $(H'_3)$ hold true for the fixed $\varepsilon> 0$. Take $s>\frac{1}{2}$. We obtain the following assertions:
\item[1.] $\sigma_{\rm ess}(H_m)=\sigma_{\rm ess}(\Delta_{T})$.
\item[2.] The eigenvalues of $H_m$, distinct from $-\frac{1}{2}$, $-\frac{1}{3}$ and $1$ are of finite multiplicity and can accumulate
only at $-\frac{1}{2}$, $-\frac{1}{3}$ and $1$.

\item[3.] The singular continuous spectrum of $H_m$ is empty.

\item[4.] Take $[a,b]$ included in $\mathbb{R}\setminus\left(\kappa(H_m)\cup\sigma_{\rm p}(H_m)\right).$ The following limit is finite:
$$\lim_{\rho\rightarrow 0^+}\sup_{\lambda\in[a,b]}\|\Lambda^{-s}(Q)(H_m-\lambda-\rmi\rho)^{-1}\Lambda^{-s}(Q)\|<\infty.$$
Moreover, in the norm topology of bounded operators, the boundary values of the resolvent:
\[ [a,b] \ni\lambda\mapsto\lim_{\rho\to0^{\pm}}\Lambda^{-s}(Q)(H_m-\lambda-\rmi\rho)^{-1}\Lambda^{-s}(Q) \mbox{ exists and is continuous}.\]

\item[5.] There exists $c>0$ such that for all $f\in\ell^2(\Vc,m)$, we have:
\[\int_{\R}\|\Lambda^{-s}(Q)e^{-\rmi tH_m}E_{[a,b]}(H_m)f\|^2dt\leq c\|f\|^2,\]
where $E_{[a,b]}(H_m)$ is the spectral projection of $H_m$ above $[a,b].$
\end{theorem}
In point $1.$,  we only need the hypotheses $(H_0)$ and $(H'_0)$.  
Points $2.$$-5.$ are standard consequences of Mourre's theory, where we establish a Mourre estimate and verify the hypotheses of regularity.
We refer to Section \ref{s:mourre}, for historical references and an introduction on the subject. Point $4.$  is called a \emph{Limiting Absorption Principle}. It implies that the spectrum is purely absolutely continuous above $\mathbb{R}\setminus\left(\kappa(H_m)\cup \sigma_{\rm p}(H_m)\right)$. 
Specifically, Riemann Lebesgue's Theorem ensures that the solution to the Schr{\"o}dinger equation escapes at infinity. That is, for $f$ belonging to the absolutely continuous subspace of $\Delta_m$ and $n\in\Vc,$ \begin{align}\label{1}\lim_{|t|\rightarrow\infty}\left(e^{\rmi t\Delta_m}f\right)(n)=0.\end{align} While \eqref{1} confirms that the particle escapes at infinity. Point $5.$ indicates that the particle concentrates where $\Lambda^{s}$ is large. Point $5.$ corresponds to the fact that $\Lambda^{s}$ is locally $H_m$-smooth over $[a,b]$, e.g. \cite[Section VIII.C]{RS}.

The concrete framework of this work allows us to explicitly define the set of the critical points $\left\{-\frac{1}{2},\ -\frac{1}{3},\ 1\right\}$, which corresponds to the energy where, after Fourier transform, the symbol of $\Delta_T$ is zero, at this energies see Lemma \ref{Fourier}. Intuitively, there is no propagation, see Lemma \ref{c}. In \cite{PR}, the authors use a general and abstract Floquet-Bloch approach, which ensures the existence of critical points via direct integral decomposition, see also \cite{GN} for a general theory. However, they do not give this set explicitly.

We now, give the structure of our paper. Section \ref{s:mourre} presents a brief overview of Mourre theory. Subsection \ref{LTL} studies the model and proves the Mourre estimate for the Laplacian acting on a triangular lattice. Subsection \ref{section7} examines metric perturbation and the addition of a potential.
Finally, Subsection \ref{P.result} establishes the main results of Theorem \ref{t:LAP}. \\
\textbf{Acknowledgements:}
We would like to thank the reviewers for their precious comments on the manuscript.   
\section{The Mourre theory}\label{s:mourre}
In 1956, C.R. Putnam provided a condition for the spectrum of a self-adjoint operator $H$ to be purely absolutely continuous, assuming the existence of a bounded self-adjoint operator $B$ such that $[H,\rmi B]>0$. However,      
the boundedness of $B$ is a strong constraint for applications. The Mourre theory has attracted significant interest
since its introduction in 1980 (cf., \cite{Mo81,Mo83}). Many works have proved the importance of the Mourre
commutator theory for the point and continuous spectra of a sufficiently broad class of self-adjoint operators.
Among the interesting works, we can see \cite{CGH,GGM1,GG,JMP,S}, the book \cite{ABG}, the master courses \cite{G} and more recent results such as
\cite{GJ,Ge,GoMa}.

Now, we recall
Mourre's commutator theory. Let $H$ and $A$ be two self-adjoint operators acting on a complex Hilbert
space $\Hc$. Suppose also $H\in\Bc(\Hc)$. We denote by $\|\cdot\|$ the norm of bounded operators on $\Hc$. Thanks to the operator $A$, we
study several spectral properties of $H$.
Given $k\in \N$,  we say that $H\in \Cc^k(A)$ if for all $f\in \Hc$ the
map $\R\ni t\mapsto e^{\rm i t A}H e^{-\rm i t A}f$ has the usual $\Cc^k(\mathbb{R})$ regularity, i.e. $\R\ni t\mapsto e^{\rm i t A}H e^{-\rm i t A}\in \Bc(\Hc)$ has the usual $\Cc^k(\mathbb{R})$ regularity with  $\Bc(\Hc)$ endowed with the strong operator topology.
We say that $H\in \Cc^{k,u}(A)$, if the map $\R\ni t\mapsto e^{\rm i t A}H e^{-\rm i t A}\in\Bc(\Hc)$ has the usual $\Cc^k(\mathbb{R})$ regularity, with $\Bc(\Hc)$ endowed with the norm operator topology.
The form $[H,\rmi A]$ is defined on $\Dc(A)\times\Dc(A)$ by $\langle f,[H,\rmi A]g\rangle:=\rmi\left(\langle Hf,Ag\rangle+\langle Af,Hg\rangle\right)$. By \cite[Lemma 6.2.9]{ABG} $H\in\Cc^1(A)$ if and only if the form $[H,\rmi A]$ extends to a bounded operator in which case we denote by $[H,\rmi A]_{\circ}$.  
We say that $H\in\Cc^{0,1}(A)$ if
\[\int^1_0\| [H,e^{\rm i t A}]\| \frac{dt}{t}<\infty\]
and that $H\in\Cc^{1,1}(A)$ if
\[\int^1_0\| [[H,e^{\rm i t A}],e^{\rm i t A}]\| \frac{dt}{t^2}<\infty.\]
Thanks to \cite[p. 205]{ABG}, we have the following of vector spaces inclusions:
\begin{align}\label{classe}\Cc^2(A)\subset\Cc^{1,1}(A)\subset\Cc^{1,u}(A)\subset\Cc^1(A)\subset\Cc^{0,1}(A).\end{align} 
Note that, for a bounded operator $H$, if $[H,\rmi A]_{\circ}\in \Cc^{0,1}(A)$ then $H\in \Cc^{1,1}(A).$

The \emph{Mourre estimate} for $H$ on an open interval $\Ic$ of $\R$ holds true if there exist $c>0$ and a compact operator $K$ such that:
\begin{align}\label{eqmourre}
E_\Ic (H)[H, \rmi A]_{\circ}E_\Ic(H)\geq E_\Ic (H)(c\, +\, K) E_\Ic (H),
\end{align}
where $E_{\Ic}(H)$ is the spectral measure of $H$ above $\Ic$.
Mourre's commutator theory aims to
prove a \emph{Limiting Absorption Principle} (LAP), see \cite[Theorem 7.6.8]{ABG}.
\begin{theorem}Let $H$ be a self-adjoint operator, with $\sigma(H)\neq \R$. Assume that $H\in\Cc^1(A)$ and the Mourre estimate \eqref{eqmourre} holds true for $H$ on $\Ic$. Then:
\item[1.] If $K=0$, then $H$ has no eigenvalues in $\Ic$.
\item[2.] The number of eigenvalues of $H$ on $\Ic$ counted with multiplicity is finite.
\item[3.] If $H\in\Cc^{1,1}(A)$, $s>1/2$ and $\Ic'$ a compact sub-interval of $\Ic$ that contains no eigenvalue, then \[\sup_{\Re(z)\in\Ic',\Im(z)\neq0}\|\langle A\rangle^{-s}(H-z)^{-1}\langle A\rangle^{-s}\| \mbox{ is finite}.\]
\item[4.] In the norm topology of bounded operators, the boundary values of the resolvent:
\[\Ic'\ni\lambda\mapsto\lim_{\rho\to0^{\pm}}\langle A\rangle^{-s}(H-\lambda-\rmi\mu)^{-1}\langle A\rangle^{-s} \mbox{ exists and is continuous}.\]
\end{theorem}
For more details, see \cite[Proposition 7.2.10, Corollary 7.2.11, Theorem 7.5.2]{ABG}.
\section{Proof of the main result}\label{section3}
We aim to prove the Theorem \ref{t:LAP}. Subsection \ref{LTL} studies the Laplacian on a triangular lattice and proves its Mourre estimate. Subsection \ref{section7}, examines the metric perturbations and addition of potential. Finally, Subsection \ref{P.result} proves Theorem \ref{t:LAP}.\\  
\subsection{Laplacian on the triangular lattice}\label{LTL}

Given $f\in\ell^2(\mathbb{Z})$, we set 
$$U_1f(n)=U_2f(n):=f(n-1).\text{ Note that } U_1^*f(n)=U_2^*f(n)=f(n+1).$$ 
Under the identification
\[\ell^2(\mathbb{Z}^2,1)\simeq \ell^2(\mathbb{Z},1)\otimes \ell^2(\mathbb{Z},1).\]
For all $f\in\ell^2(\mathbb{Z}^2,1)$, we have: 
\begin{align*}
&\left(U_1\otimes 1\right) f(n_1,n_2)=f(n_1-1,n_2),\ \left(U_1\otimes 1\right)^{*} f(n_1,n_2)=f(n_1+1,n_2),\\
&\left(1\otimes U_2\right)  f(n_1,n_2)=f(n_1,n_2-1),\ \left(1\otimes U_2\right)^{*} f(n_1,n_2)=f(n_1,n_2+1),\\
&\left(U_1\otimes U_2^*\right) f(n_1,n_2)=f(n_1-1,n_2+1),\ \left(U_1^{*}\otimes U_2\right) f(n_1,n_2)=f(n_1+1,n_2-1).
\end{align*}
Note that $(U_1\otimes 1)^*=U_1^*\otimes 1$ and $(1\otimes U_2)^*=1\otimes U_2^*.$

Let $\Sc:=\left\{f:\mathbb{Z}^2\rightarrow \mathbb{C}\mbox{ such that,\ for all } N\in\mathbb{N}\ \sup_{n}\left|(1+n_1^2+n_2^2)^Nf(n)\right|<\infty\right\}$, it is the discrete Schwartz space.
For all $f\in \Sc$, we have: $$[U_1^*\otimes 1,Q_1\otimes 1]f(n)=U^*_1\otimes 1f(n),$$ then by density we have:  
\begin{equation}\label{U1*}
[U_1^*\otimes 1,Q_1\otimes 1]_{\circ}f(n)=U^*_1\otimes 1f(n), \forall f\in \ell^2(\Vc,1).
\end{equation}
In the same way, we have: 
\begin{equation}\label{U1}
[U_1\otimes 1,Q_1\otimes 1]_{\circ}f(n)=-(U_1\otimes 1)f(n),
\end{equation}
\begin{equation}\label{U2}
[1\otimes U_2^*,1\otimes Q_2]_{\circ}f(n)=1\otimes U^*_2f(n),
\end{equation}
\begin{equation}\label{U2*}
[1\otimes U_2,1\otimes Q_2]_{\circ}f(n)=-(1\otimes U_2)f(n),
\end{equation}
\begin{equation}\label{U1U2*}
[1\otimes Q_2,U_1\otimes U_2^*]_{\circ}f(n)=-(U_1\otimes U_2^*)f(n),
\end{equation}
\begin{equation}\label{U1U2*Q1}
[1\otimes Q_1,U_1\otimes U_2^*]_{\circ}f(n)=U_1\otimes U_2^*f(n),
\end{equation}
\begin{equation}\label{U1*U2}
[1\otimes Q_2,U_1^*\otimes U_2]_{\circ}f(n)=U_1^*\otimes U_2f(n),
\end{equation}
and
\begin{equation}\label{U1*U2Q1}
[1\otimes Q_1,U_1^*\otimes U_2]_{\circ}f(n)=-(U_1^*\otimes U_2)f(n).
\end{equation}
We denote by $\Cc^{\infty}_{2\pi}(|-\pi,\pi]^2)$ the set of functions defined on $[-\pi,\pi]^2$ that are of class $\Cc^{\infty}$ and $2\pi$-periodic.

First, we rewrite the Laplacian on a triangular lattice.
\begin{lemma}\label{Laplacian}
The Laplacian on a triangular lattice is given by:
\begin{align*}
\Delta_{T}
:=&\frac{1}{6}\Bigg(U_1\otimes 1+U_1^{*}\otimes 1+1\otimes U_2+1\otimes U_2^{*}+U^*_1\otimes U_2
+U_1\otimes U^{*}_2\Bigg).
\end{align*}
\end{lemma}
\proof
Recalling \eqref{Na} and \eqref{map}. Let $f\in \ell^2(\mathbb{Z}^2,1)$, we have:
\begin{align*}
(\Delta_{T}f)(n)
=&\frac{1}{6}\Big(f(n+v_1)+f(n-v_1)+f(n+v_2)+f(n-v_2)\\
&+f(n+(v_1-v_2))+f(n-(v_1-v_2))\Big)\\
=&\frac{1}{6}\Big(f(n_1+1,n_2)+f(n_1-1,n_2)+f(n_1,n_2+1)\\
&+f(n_1,n_2-1)+f(n_1+1,n_2-1)+f(n_1-1,n_2+1)\Big).
\end{align*}
This gives the result.
\qed

Now, we define the Fourier transform $\Fc:\ell^2(\mathbb{Z}^2,1)\longrightarrow L^2([-\pi,\pi]^2)$ through
\[\Fc f(x):=\frac{1}{2\pi}\sum_n f(n)e^{-\rmi\langle n,x\rangle},\ \forall f\in\ell^2(\mathbb{Z}^2,1). \]
The inverse Fourier transform $\Fc^{-1} :L^2([-\pi,\pi]^2)\longrightarrow\ell^2(\mathbb{Z}^2,1)$ is given by
\[\Fc^{-1} f(n)=\frac{1}{2\pi}\int_{[-\pi,\pi]^2} f(x)e^{\rmi\langle n,x\rangle}dx,\ \forall f\in L^2([-\pi,\pi]^2).\]

\begin{lemma}\label{Fourier}
For $f\in L^2([-\pi,\pi]^2)$, we have:
\begin{align*}\Fc\Delta_{T} \Fc^{-1}f(x):=(F(Q)f)(x)=F(x)f(x),\end{align*}
with
\begin{align*}
F(x):=\frac{1}{3}\left(\cos(x_1)+\cos(x_2)+\cos(x_1-x_2)\right),
\end{align*}
where $x:=(x_1,x_2).$
\end{lemma}
\proof
Let $f\in L^2([-\pi,\pi]^2)$, we have:
\begin{align*}
\Fc(U_1\otimes1\ \Fc^{-1}f)(x)
=&\frac{1}{2\pi}\sum_{n}(U_1\otimes 1\ \Fc^{-1}f)(n)e^{-\rmi\langle x, n\rangle}\\
=&\frac{1}{2\pi}\sum_{n}(\Fc^{-1}f)(n_1-1,n_2)e^{-\rmi\langle x, n\rangle}\\
=&\frac{1}{2\pi}\sum_{n}(\Fc^{-1}f)(n_1,n_2)e^{-\rmi\left(\langle x_1, n_1+1\rangle+\langle x_2,n_2\rangle\right)}\\
=&\frac{1}{2\pi}\sum_{n}(\Fc^{-1}f)(n_1,n_2)e^{-\rmi x_1}e^{-\rmi\langle x, n\rangle}.
\end{align*}
Then,
$\Fc (U_1\otimes1\ \Fc^{-1}f)(x)=e^{-\rmi x_1}f(x)$. The other terms are treated in the same way. We obtain the result.
\qed

Next, we compute to the spectrum.
\begin{lemma} \label{specter}
$\sigma(F(Q))=\left[-\frac{1}{2},1\right].$
\end{lemma}
\proof 
From the definition of $F$, we see that $F \leq 1 = F(0,0)$.\newline
On $[-\pi, \pi]^2$, we introduce the functions:
\begin{equation*}
    X = \cos\left( \frac{x_1 + x_2}{2} \right), \quad Y = \cos\left( \frac{x_1 - x_2}{2} \right).
\end{equation*}
Using trigonometric formulas, we obtain:
\begin{equation*}
    3F = 2XY + 2Y^2 - 1.
\end{equation*}
We observe that:
\begin{equation}\label{13}
    3F + \frac{3}{2} = 2\left(Y + \frac{X}{2}\right)^2 + \frac{1 - X^2}{2} \geq 0,
\end{equation}
since $|X| \leq 1$. This shows that $F \geq -\frac{1}{2}$. Thus, the range of $F$ is included in $[-\frac{1}{2},1]$.
Now, using \eqref{13}, we have:
\begin{equation}\label{14}
    F = -\frac{1}{2} \Longleftrightarrow |X| = 1 \text{ and } Y = -\frac{X}{2}.
\end{equation}
If $X = -1$, then $x_1 = x_2 = \pi$ or $x_1 = x_2 = -\pi$ thus $Y = 1$ (contradiction with \eqref{14}). Therefore:
\begin{equation*}
    F = -\frac{1}{2} \Longleftrightarrow X = 1 \text{ and } Y = -\frac{1}{2}.
\end{equation*}
Taking into account that $(x_1, x_2)$ lives in $[-\pi, \pi]^2$, we obtain:
\begin{equation*}
    F = -\frac{1}{2} \Longleftrightarrow (x_1, x_2) \in \left\{\left(-\frac{2\pi}{3}, \frac{2\pi}{3}\right), \left(\frac{2\pi}{3}, -\frac{2\pi}{3}\right)\right\}.
\end{equation*}
Thus $1$ and $-\frac{1}{2}$ belongs to the range of $F$ and, since $F$ is continuous, the range of $F$ is
$[-\frac{1}{2}, 1]$. This gives the spectrum of $F(Q)$ (e.g. \cite[Vol 1, p. 229]{RS}).\qed
\begin{lemma}
The critical values of $F$ are $-\frac{1}{2}$, $-\frac{1}{3}$, and $1$.
\end{lemma} 
\proof
$F$ is a smooth. We have:
\begin{equation*}
    3\nabla F = 2Y \nabla X + 2(X + 2Y) \nabla Y,
\end{equation*}
where $X$ and $Y$ are defined in the proof of Lemma \ref{specter} and  
\begin{equation*}
    \nabla X = -\frac{1}{2} \sin\left(\frac{x_1 + x_2}{2}\right) \left(
                                                                   \begin{array}{c}
                                                                     1 \\
                                                                     1 \\
                                                                   \end{array}
                                                                 \right)
    , \quad
    \nabla Y = \frac{1}{2} \sin\left(\frac{x_1 - x_2}{2}\right)\left(
                                                                 \begin{array}{c}
                                                                   -1 \\
                                                                   1 \\
                                                                 \end{array}
                                                               \right)
     .
\end{equation*}
In particular, $\nabla X$ and $\nabla Y$ are orthogonal. Thus, they are independent vectors unless one of them is zero. Therefore:
\begin{align*}
    &\nabla F = 0 \Longleftrightarrow 
        \left(\nabla X = 0 = \nabla Y\right) 
        \text{or } \left(\nabla X = 0, \nabla Y \neq 0, \text{ and } X + 2Y = 0\right) \\
        &\hspace*{2,1cm}\text{or } \left(\nabla Y = 0, \nabla X \neq 0, \text{ and } Y = 0\right) \\
        &\hspace*{2,1cm}\text{or } \left(\nabla X, \nabla Y \text{ are independent, } Y = 0, X + 2Y = 0\right).
\end{align*}
We observe that we cannot simultaneously have $Y = 0$ and $\nabla Y = 0$.\newline
In the case of independent $\nabla X$ and $\nabla Y$, we have $X = Y = 0$, hence $F = -\frac{1}{3}$. This actually occurs for $(x_1, x_2) = (0, \pi)$.\newline
In the case $\nabla X = 0 = \nabla Y$, we must have $|X| = |Y| = 1$. Then $F = 1$ if $XY > 0$ and $F = -\frac{1}{3}$ if $XY < 0$. The first case occurs at $(0,0)$.\newline
In the last case, we have $\nabla X = 0$, $\nabla Y \neq 0$, and $X + 2Y = 0$. Since $\nabla Y \neq 0$, $x_1 \neq x_2$ thus $-\pi < \frac{(x_1 + x_2)}{2} < \pi$, yielding $X > -1$. Since $\nabla X = 0$, $|X| = 1$. Thus, $X = 1$ and $Y = -\frac{1}{2}$. Then $F = -\frac{1}{2}$, by \eqref{14}. This case occurs for $(x_1, x_2) = (-\frac{2\pi}{3}, \frac{2\pi}{3})$.\newline
We have shown that the critical values of $F$ are $-\frac{1}{2}$, $-\frac{1}{3}$, and $1$.\qed\\

We define now the conjugate operator $A$ on the discrete Schwartz space $\Sc.$
\begin{definition}
On $\Sc $, we set
\begin{align}\label{a}
A
\nonumber:=&\frac{\rmi}{6}\Bigg(\left(\frac{U_1^*-U_1}{2}\otimes 1+\frac{U_1^*\otimes U_2}{2}-\frac{U_1\otimes U_2^*}{2}\right)Q_1\\
&+\left(1\otimes \frac{U_2^*-U_2}{2}+\frac{U_1\otimes U_2^*}{2}-\frac{U_1^*\otimes U_2}{2}\right)Q_2\Bigg)\\
\nonumber&+\rmadj.
\end{align}
\end{definition}
On $\Sc$, we can rewrite $A$ as follows: 
\begin{align}\label{A}
\nonumber Af(n)
\nonumber=&\frac{\rmi}{6}\Big((Q_1+\frac{1}{2})U_1^*\otimes1-(Q_1-\frac{1}{2})U_1\otimes1\\
&+1\otimes(Q_2+\frac{1}{2})U_2^*-1\otimes(Q_2-\frac{1}{2})U_2\\
\nonumber&+(Q_2-Q_1+1)U_1\otimes U_2^*+(Q_1-Q_2+1)U_1^*\otimes U_2\Big)f(n).
\end{align}
\begin{lemma}\label{Alambda}
There exists $C> 0$ such that for all $f\in\Sc$, we have: \[\left\|Af\right\|^2\leq C\left\|\Lambda(Q) f\right\|^2.\]
In particular, there is $C'>0$ such that $$\|\langle A\rangle f\|^2\leq C'\|\Lambda(Q)f\|^2,\ \forall f\in\Sc.$$ 
\end{lemma}
\proof
Let $f\in\Sc$. Here all constants are denoted by $C$ and are independent of $f$. 
\begin{align*}
\|Af\|^2
=&\sum_{n\in\Vc}\Big|\frac{\rmi}{6}\Big((Q_1+\frac{1}{2})U_1^*\otimes1-(Q_1-\frac{1}{2})U_1\otimes1\\
&+1\otimes(Q_2+\frac{1}{2})U_2^*-1\otimes(Q_2-\frac{1}{2})U_2\\
&+(Q_2-Q_1+1)U_1\otimes U_2^*+(Q_1-Q_2+1)U_1^*\otimes U_2\Big)f(n)\Big|^2\\
\leq& C\sum_n\Big|(Q_1+\frac{1}{2})U_1^*\otimes1f(n)\Big|^2+\Big|(Q_1-\frac{1}{2})U_1\otimes1f(n)\Big|^2\\
&+\Big|1\otimes(Q_2+\frac{1}{2})U_2^*f(n)\Big|^2+\Big|1\otimes(Q_2-\frac{1}{2})U_2f(n)\Big|^2\\
&+\Big|(Q_2-Q_1+1)U_1\otimes U_2^*f(n)\Big|^2+\Big|(Q_1-Q_2+1)U_1^*\otimes U_2f(n)\Big|^2.
\end{align*}
We treat the first term of $\|Af\|^2$, we have:
\begin{align*}
&\sum_{n}\left|(Q_1+\frac{1}{2})U_1^*\otimes1f(n)\right|^2\\
=&\sum_{n}\left|\left(U_1^*\otimes1(Q_1+\frac{1}{2})+\left[Q_1+\frac{1}{2},U_1^*\otimes1\right]_{\circ}\right)f(n)\right|^2,\text{ by } \eqref{U1*}\\
\leq& C\left(\|(Q_1+\frac{1}{2})f\|^2+\|f\|^2\right)
=C\left(\langle f,(Q_1+\frac{1}{2})^2f\rangle+\|f\|^2\right)\\
\leq& C\left(\|(Q_1^2+\frac{1}{2})^{\frac{1}{2}}f\|^2+\|f\|^2\right)
\leq C\|\Lambda(Q) f\|^2
\end{align*}
and we estimate the next term
\begin{align*}
&\sum_{n\in\Vc}|(Q_2-Q_1+1)U_1\otimes U_2^*f(n)|^2\\
=&\sum_n\left|\left((Q_2+\frac{1}{2})-(Q_1+\frac{1}{2})+1\right)U_1\otimes U_2^*f(n)\right|^2\\
=&\sum_{n}\Big|\Big(U_1\otimes U_2^*(Q_2+\frac{1}{2})+\left[Q_2+\frac{1}{2},U_1\otimes U_2^*\right]_{\circ}+U_1\otimes U_2^*(Q_1+\frac{1}{2})\\
&+\left[Q_1+\frac{1}{2},U_1\otimes U_2^*\right]_{\circ}+U_1\otimes U_2^*
\Big)f(n)\Big|^2,\text{ by } \eqref{U1U2*} \text{ and }\eqref{U1U2*Q1}\\
\leq& C\sum_n |(Q_2+\frac{1}{2})f(n)|^2+|(Q_1+\frac{1}{2})f(n)|^2+|f(n)|^2\\
=&C\left(\langle f,(Q_2+\frac{1}{2})^2f\rangle+\langle f,(Q_1+\frac{1}{2})^2f\rangle+\|f\|^2\right)\\
\leq& C\left(\langle f,(Q_2^2+\frac{1}{2})f\rangle+\langle f,(Q_1^2+\frac{1}{2})f\rangle+\|f\|^2\right)\\
\leq& C\left(\|(Q_2^2+\frac{1}{2})^{\frac{1}{2}}f\|^2+\|(Q_1^2+\frac{1}{2})^{\frac{1}{2}}f\|^2+\|f\|^2\right)
\leq C\|\Lambda(Q) f\|^2.
\end{align*}
The rest of the terms are bounded in the same way by using \eqref{U1}, \eqref{U2}, \eqref{U2*}, \eqref{U1*U2} and \eqref{U1*U2Q1}.
This gives the first point. Next, given $f\in\Sc$ note that \begin{align*}\|\langle A\rangle f\|^2=\langle\langle A\rangle f, \langle A \rangle f\rangle=
\langle f,(1+A^2)f\rangle=\|f\|^2+\|Af\|^2.\end{align*}This concludes the proof.
\qed
\begin{remark}\label{Lambdavarepsilon}
Thanks to Lemma \ref{Alambda} and since $\|\langle A\rangle^0f\|^2\leq \|\Lambda^0 (Q)f\|^2, \text{ for all } f\in\Sc$, by real interpolation, e.g. \cite[Theorem 4.1.2,p.88]{BL}, for all $\gamma\in [0,1]$ there is $C_{\gamma}$ such that \[\left\|\langle A\rangle^{\gamma}f\right\|^2\leq C_{\gamma}\left\|\Lambda^{\gamma}(Q) f\right\|^2, \text{ for all } f\in\Sc.\]  
\end{remark}
\begin{lemma}\label{Nelson}
$A$ is essentially self-adjoint on $\Sc$. We keep the notation $A$ for its closure in the sequel.
\end{lemma}
\proof
First, by definition, see \eqref{a}, $A$ is symmetric operator on $\Sc$. By Lemma \ref{Alambda}, there exists $C$ such that for all $f\in\Sc$, we have: \[\left\|Af\right\|^2\leq C\left\|\Lambda(Q) f\right\|^2.\] By the Nelson's Lemma, e.g. \cite[Theorem X.37]{RS}, it suffices to prove 
 \[\exists C>0,\ \forall f\in\Sc,\ \Big|\langle f,[\Lambda(Q),A]f\rangle\Big|\leq C\big\|\Lambda^\frac{1}{2}(Q)f\big\|.\] to ensure to that $A$, defined on $\Sc$, extends to a self-adjoint operator.  
 
Let $f\in\Sc$. We denote all constants by $C$, we infer:
\begin{align*}
[\Lambda(Q),A]f(x)
=&-\frac{\rmi}{6}\Big((Q_1+\frac{1}{2})L_1(Q_1)U_1^*\otimes 1
-(Q_1-\frac{1}{2})L_2(Q_1)U_1\otimes 1\\
&+(Q_2+\frac{1}{2})L_3(Q_2)1\otimes U_2^*
-(Q_2-\frac{1}{2})L_4(Q_2)1\otimes U_2\\
&+(Q_2-Q_1+1)L_5(Q_1,Q_2)U_1\otimes U_2^*\\
&+(Q_1-Q_2+1)L_6(Q_1,Q_2)U_1^*\otimes U_2\Big)f(x),
\end{align*}
with
\begin{align*}
L_1(Q_1):=(Q_1^2+\frac{1}{2})^{\frac{1}{2}}-(Q_1^2+2Q_1+\frac{3}{2})^{\frac{1}{2}},\end{align*}
\begin{align*} L_2(Q_1):=(Q_1^2+\frac{1}{2})^{\frac{1}{2}}-(Q_1^2-2Q_1+\frac{3}{2})^{\frac{1}{2}},\end{align*}
\begin{align*} L_3(Q_2):=(Q_2^2+\frac{1}{2})^{\frac{1}{2}}-(Q_2^2+2Q_2+\frac{3}{2})^{\frac{1}{2}},\end{align*}
\begin{align*} L_4(Q_2):=(Q_2^2+\frac{1}{2})^{\frac{1}{2}}-(Q_2^2-2Q_2+\frac{3}{2})^{\frac{1}{2}},\end{align*}
\begin{align*} L_5(Q_1,Q_2)
:=&(Q_1^2+\frac{1}{2})^{\frac{1}{2}}+(Q_2^2+\frac{1}{2})^{\frac{1}{2}}-(Q_1^2-2Q_1+\frac{3}{2})^{\frac{1}{2}}\\
&-(Q_2^2+2Q_2+\frac{3}{2})^{\frac{1}{2}}\end{align*}
and
\begin{align*} L_6(Q_1,Q_2)
:=&(Q_1^2+\frac{1}{2})^{\frac{1}{2}}+(Q_2^2+\frac{1}{2})^{\frac{1}{2}}-(Q_1^2+2Q_1+\frac{3}{2})^{\frac{1}{2}}\\
&-(Q_2^2-2Q_2+\frac{3}{2})^{\frac{1}{2}}.\end{align*}
We estimate the first term of $\Big|\langle f,[\Lambda(Q),A]f\rangle\Big|$, we have:
\begin{align*}
&\left|\sum_{n\in\Vc}\overline{f(n)}(n_1+\frac{1}{2})L_1(n_1)U_1^*\otimes 1f(n)\right|\\ \leq&\sum_n\left|n_1+\frac{1}{2}\right|^{\frac{1}{2}}\cdot\left|f(n)\right|\cdot\left|L_1(n_1)\right|\cdot\left|n_1+\frac{1}{2}\right|^{\frac{1}{2}}\cdot
\left|U_1^*\otimes 1f(n)\right|\\
\leq& \left\|\left|Q_1+\frac{1}{2}\right|^{\frac{1}{2}}f\right\|\cdot\left\|\frac{Q_1+1}{(Q_1^2+\frac{1}{2})^{\frac{1}{2}}+(Q_1^2+2Q_1+\frac{3}{2})^{\frac{1}{2}}}
\left|Q_1+\frac{1}{2}\right|^{\frac{1}{2}}U_1^*\otimes 1f(Q_1,Q_2)\right\|\\
\leq& C\left\|\left|Q_1+\frac{1}{2}\right|^{\frac{1}{2}}f\right\|\cdot\left\|\left|Q_1+\frac{1}{2}\right|^{\frac{1}{2}}U_1^*\otimes 1f(Q_1,Q_2)\right\|\\
\leq& C\left(\left\|\left|Q_1+\frac{1}{2}\right|^{\frac{1}{2}}f\right\|^2+\left\|\left|Q_1+\frac{1}{2}\right|^{\frac{1}{2}}U_1^*\otimes 1f(Q_1,Q_2)\right\|^2\right)\\
\leq& C\left(\left\|\left|Q_1+\frac{1}{2}\right|^{\frac{1}{2}}f\right\|^2+\left\|f\right\|^2\right)\\
\leq& C\left(\left\|\left(Q_1^2+\frac{1}{2}\right)^{\frac{1}{4}}f\right\|^2+\left\|f\right\|^2\right)\leq C\left\|\Lambda^{\frac{1}{2}}(Q)f\right\|^2
\end{align*}
and we treat the last term
\begin{align*}
&\left|\sum_{n\in\Vc}\overline{f(n)}\left(n_1-n_2+1\right)L_5(n_1,n_2)U_1^*\otimes U_2 f(n)\right|\\
=&\sum_n\left|\overline{f(n)}\left(\left(n_1+\frac{1}{2}\right)+\left(n_2+\frac{1}{2}\right)+1\right)L_5(n_1,n_2)U_1^*\otimes U_2 f(n)\right|\\
\leq& \sum_n\left|\overline{f(n)}\left(n_1+\frac{1}{2}\right)L_5(n_1,n_2)U_1^*\otimes U_2 f(n)\right|\\
&+\sum_n\left|\overline{f(n)}\left(n_2+\frac{1}{2}\right)L_5(n_1,n_2)U_1^*\otimes U_2 f(n)\right|\\
&+\sum_n\left|\overline{f(n)}L_5(n_1,n_2)U_1^*\otimes U_2 f(n)\right|\\
\leq& C\left(\left\|\left|Q_1+\frac{1}{2}\right|^{\frac{1}{2}}f\right\|^2+\left\|\left|Q_2+\frac{1}{2}\right|^{\frac{1}{2}}f\right\|^2+\left\|f\right\|^2\right)\\
\leq& C\left(\left\|\left(Q_1^2+\frac{1}{2}\right)^{\frac{1}{4}}f\right\|^2+\left\|\left(Q_2^2+\frac{1}{2}\right)^{\frac{1}{4}}f\right\|^2+
\left\|f\right\|^2\right)\leq C\left\|\Lambda^{\frac{1}{2}}(Q)f\right\|^2.
\end{align*}
The other terms are controlled in the same way. This gives:
\begin{align*}
\left|\left\langle f,[\Lambda(Q),A]f\right\rangle\right|\leq C\left\|\Lambda^{\frac{1}{2}}(Q)f\right\|^2.
\end{align*}
As $\Sc$ is a core for $\Lambda(Q)$, applying \cite[Theorem X.37]{RS}, the result follows.
\qed\\

Let $f\in \Cc_{2\pi}^{\infty}([-\pi,\pi]^2)$, we set:
\begin{align*}
\widehat{A}f:
&=\frac{\rmi}{2}\left(\nabla F(Q_1,Q_2)\cdot\nabla+\nabla\cdot\nabla F(Q_1,Q_2)\right)f.
\end{align*}
\begin{lemma}\label{Operator}
On $\Cc_{2\pi}^{\infty}([-\pi,\pi]^2)$, we have:
\begin{align*}
\widehat{A}
=&\frac{\rmi}{6}\Bigg(\left(-\sin(Q_1)-\sin(Q_1-Q_2)\right)\frac{\partial}{\partial x_1}\\
&+
\left(-\sin(Q_2)+\sin(Q_1-Q_2)\right)
\frac{\partial}{\partial x_2}
\Bigg)+\rmadj.
\end{align*}
\end{lemma}
\proof
Let $f\in \Cc_{2\pi}^{\infty}([-\pi,\pi]^2)$, we have:
\[\begin{pmatrix}
                               \frac{\partial F}{\partial x_1}\\
                               \frac{\partial F}{\partial x_2}\\
                            \end{pmatrix}(x_1,x_2)=\frac{1}{3}\begin{pmatrix}
                              -\sin(x_1)-\sin(x_1-x_2) \\
                              -\sin(x_2)+\sin(x_1-x_2) \\
                            \end{pmatrix}\]
and
\begin{align*}
\left\langle\begin{pmatrix}
                               \frac{\partial F}{\partial x_1}\\
                               \frac{\partial F}{\partial x_2} \\
                            \end{pmatrix},\left(
                     \begin{array}{c}
                       \frac{\partial f}{\partial x_1} \\
                       \frac{\partial f}{\partial x_2} \\
                     \end{array}
                  \right)\right\rangle(x_1,x_2)
=&\frac{1}{3}(\left(-\sin(x_1)-\sin(x_1-x_2)\right)\frac{\partial f}{\partial x_1}(x_1,x_2)\\
&+\left(-\sin(x_2)+\sin(x_1-x_2)\right)
\frac{\partial f}{\partial x_2}(x_1,x_2)).\end{align*}
This concludes the result.
\qed
\begin{lemma}\label{regilarite}
On $\Cc^{\infty}_{2\pi}([-\pi,\pi]^2)$,
we have:
\begin{align}\label{F}
\widehat{A}
=\Fc A \Fc^{-1}.
\end{align}
\end{lemma}
\proof
We recall \eqref{a}. Let $f\in\Cc^{\infty}_{2\pi}([-\pi,\pi]^2)$, 
we infer:
\begin{align*}
\frac{1}{2\rmi}\left( \Fc\left((U_1^*-U_1)\otimes 1\right)\Fc^{-1} f\right)(x)
=&\frac{1}{2\rmi}\left(e^{\rmi x_1}-e^{-\rmi x_1}\right)f(x)=\sin(x_1)f(x),
\end{align*}
\[\frac{1}{2\rmi}\left(\Fc\left(U_1\otimes U_2^*-U_1^*\otimes U_2\right)\Fc^{-1}f\right)(x)=\sin(x_1-x_2)f(x)\]
and
\begin{align*}
\left(\Fc(-\rmi Q_1)\Fc^{-1} f\right)(x)
=\frac{\partial f}{\partial x_1}(x).
\end{align*}
The other terms are estimated similarly, we obtain the result.
\qed
\begin{remark}
Since $\Fc(\Sc)=\Cc_{2\pi}^{\infty}\left([-\pi,\pi]^2\right)$ and recall Lemma \ref{regilarite}, by density we have $\widehat{A}$ is essentially self-adjoint on $\Cc_{2\pi}^{\infty}([-\pi,\pi]^2)$ and we denote by  $\widehat{A}$ its closure. Note that \eqref{F} extends to the closure and $\Dc(A)=\Fc^{-1}\Dc(\widehat{A}).$
\end{remark}
\begin{lemma}\label{c}
We have $F(Q)\in \Cc^1(\widehat{A})$ and therefore $\Delta_{T}\in\Cc^1(A)$. Moreover,
\begin{align}\label{*}
\left[F(Q_1,Q_2),\rmi \widehat{A}\right]_{\circ}
=\Big\|\nabla F(Q_1,Q_2)\Big\|^2_{\mathbb{C}^2},
\end{align}
where $$\Big\|\nabla F(Q_1,Q_2)\Big\|^2_{\mathbb{C}^2}f(x_1,x_2)=\sum_{j=1}^{2}\left(\frac{\partial F}{\partial x_j}\right)^2(x_1,x_2)f(x_1,x_2), \forall f\in \Cc_{2\pi}^{\infty}([-\pi,\pi]^2).$$
\end{lemma}
 \proof
For $f\in \Cc_{2\pi}^{\infty}([-\pi,\pi]^2)$, we have:
\begin{align*}
\left[F(Q_1,Q_2),\rmi \widehat{A}\right]f(x_1,x_2)
=&-\left[F(Q_1,Q_2),\nabla F(Q_1,Q_2)\cdot\nabla\right]f(x_1,x_2)\\
=&-\left[F(Q),\sum_{j=1}^{2}\frac{\partial F}{\partial x_j}(Q)\frac{\partial}{\partial x_j}\right]f(x_1,x_2)\\
=&-\sum_{j=1}^{2}F(x_1,x_2)\frac{\partial F}{\partial x_j}(x_1,x_2)\frac{\partial f}{\partial x_j}(x_1,x_2)\\
&+\frac{\partial F}{\partial x_j}(x_1,x_2)\frac{\partial (Ff)}{\partial x_j}(x_1,x_2)\\
=&\sum_{j=1}^{2}\left(\frac{\partial F}{\partial x_j}\right)^2(x_1,x_2)f(x_1,x_2).
\end{align*}
As $\nabla F\in L^{\infty}([-\pi,\pi]^2), \hbox{ there exists }c> 0 \hbox{ such that } \left\|\left[F(Q_1,Q_2),\rmi \widehat{A}\right]f\right\|_{\mathbb{C}^2}\leq c\|f\|$. By density and thanks to \cite[Lemma 6.2.9]{ABG}, we obtain $F(Q_1,Q_2)\in \Cc^1(\widehat{A})$. Recalling $\Fc(\Sc)=\Cc_{2\pi}^{\infty}\left([-\pi,\pi]^2\right)$, as the Fourier transform is unitary, we obtain $\left\|\left[\Delta_T,\rmi A\right]g\right\|_{\mathbb{C}^2}\leq c\|g\|$ for all $g\in\Sc$ which ensures that $\Delta_T\in\Cc^1(A).$
\qed\\

We establish to the Mourre estimate for the unperturbed Laplacian.
\begin{proposition}\label{Mourre} 
We have $F(Q)\in \Cc^1(A)$. Moreover, let $\Ic$ be an open interval such that its closure is included in $[-\frac{1}{2},1]\backslash\{1,-\frac{1}{2},-\frac{1}{3}\}$,
there exists $c>0$ such that:
\begin{align}\label{eq:mourre}
E_\Ic (F(Q))[F(Q), \rm i \widehat{A}]_{\circ}E_\Ic(F(Q))\geq c E_\Ic(F(Q)).
\end{align}\\
Equivalently, we have: 
 \begin{align}\label{eq:mour}
E_\Ic (\Delta_{T})[\Delta_{T}, \rm i A]_{\circ}E_\Ic(\Delta_{T})\geq c E_\Ic(\Delta_{T}),
\end{align}
\end{proposition}
\proof
We work in $L^2([-\pi,\pi]^2).$ The $\Cc^1$ property is given in Lemma \ref{c}.
Let $\Ic$ be an open interval such that its closure is included in $[-\frac{1}{2},1]\backslash\{1,-\frac{1}{2},-\frac{1}{3}\}$, since $\Ic$ is bounded then by the Bolzano-Weierstrass Theorem, its closure is compact. There exists $c>0$, such that for all $(x_1,x_2)\in F^{-1}(\Ic)$, we have
$\|\nabla F(Q)\|^2_{\mathbb{C}^2}(x_1,x_2)\geq c.$ Recalling $\left[F(Q),\rmi \widehat{A}\right]_{\circ}
=\Big\|\nabla F(Q)\Big\|^2_{\mathbb{C}^2}$, by functional calculus, we have: 
$$E_\Ic (F(Q))[F(Q), \rm i \widehat{A}]_{\circ}E_\Ic(F(Q))\geq c E_\Ic(F(Q)).$$
This gives \eqref{eq:mour}, by going back to $\ell^2(\mathbb{Z}^2,1).$ \qed
\begin{proposition}\label{C2}
We have $F(Q)\in\Cc^2(\widehat{A})$ and therefore $\Delta_{T}\in\Cc^2(A).$
\end{proposition}
\proof
By Lemma \ref{c} and for $f\in \Cc^{\infty}_{2\pi}([-\pi,\pi]^2)$, we have:
\begin{align*}
\left[\left[F(Q_1,Q_2),\rmi\widehat{A}\right]_{\circ},\rmi \widehat{A}\right]f
=\left[\left\|\nabla F(Q)\right\|^2_{\mathbb{C}^2},\rmi\widehat{A}\right]f(x_1,x_2).
\end{align*}
Then,
\begin{align*}
&\left[\left\|\nabla F(Q)\right\|^2_{\mathbb{C}^2},\rmi\frac{\rmi}{6}\left(-\sin(Q_1)-\sin(Q_1-Q_2)\right)\frac{\partial }{\partial x_1}\right]f(x_1,x_2)\\
=&\frac{1}{6}\left(-\sin(x_1)-\sin(x_1-x_2)\right)\frac{\partial\left\|\nabla F\right\|^2_{\mathbb{C}^2}}{\partial x_1}(x_1,x_2)f(x_1,x_2).
\end{align*}
As $\frac{\partial\left\|\nabla F\right\|^2_{\mathbb{C}^2}}{\partial x_1}\in L^{\infty}([-\pi,\pi]^2)$,
we have:
\begin{align*}
\left\|\left[\left\|\nabla F(Q)\right\|^2_{\mathbb{C}^2},\rmi\frac{\rmi}{6}\left(-\sin(Q_1)-\sin(Q_1-Q_2)\right)\frac{\partial }{\partial x_1}\right]f\right\|^2\leq \left(\frac{10}{27}\right)^2\left\|f\right\|^2.
\end{align*}
The other terms have the same treatment. By density and thanks to \cite[Proposition 5.2.2]{ABG}, we obtain that $F(Q)\in\Cc^2(\widehat{A})$. 
As in Lemma \ref{c}, we also obtain that $\Delta_T\in\Cc^2(A).$
\qed
\subsection{The perturbed model}\label{section7}
In this subsection, we perturb the previous case by modifying the metric and adding a potential. We ask them to be, in some sense, small at infinity. We shall need some 
technicalities and start with properties of $\Lambda.$
\begin{proposition}\label{p::hypothesis}
$\Lambda$ satisfies the following assertions:
\item[1.]$\Dc(\Lambda(Q))\subset \Dc(A)$.
\item[2.]There is $c > 0$ such that for all $r > 0,\ -\rmi r$ belongs to the resolvent set of
$\Lambda$ and $r\|(\Lambda + \rmi r)^{-1}\|_{\Bc(\ell^2(\Vc,1))}\leq c.$
\item[3.]$t\rightarrow e^{\rmi t\Lambda}$ has a polynomial growth in $\ell^2(\Vc,1)$.
\item[4.] Given $\xi\in\Cc^{\infty}(\mathbb{R}.\mathbb{R})$ such that $\xi(x)=0$ near $0$ and $1$ near infinity and
$T\in\Bc(\ell^2(\Vc,1))$ symmetric, if 
\begin{align}\label{integral}
 \int_{1}^{\infty}\left\|\xi\left(\frac{\Lambda}{r}\right)T \right\|_{\Bc(\ell^2(\Vc,1))}\frac{dr}{r}<\infty
\end{align}
then $T\in \Cc^{0,1}(A).$
\end{proposition}
\proof 
\item[1.] Let $f\in\Sc$,
by Lemma \ref{Alambda} we have $\|Af\|^2\leq C\|\Lambda f\|^2$.
Since $\Lambda$ is essentially self-adjoint on $\Sc$. The result follows.
\item[2.]Note that $\Lambda$ is self-adjoint in $\ell^2(\Vc,1)$, by functional calculus it is clear, e.g \cite[Theorem VIII.5]{RS}. 
\item[3.]Again, since $\Lambda$ is self-adjoint, the norm of $t\rightarrow e^{\rmi t\Lambda}$
is $1$ for all $t\in \mathbb{R}$.
It has in particular polynomial growth.
\item[4.] Apply \cite[Proposition 7.5.7]{ABG}.\qed
\begin{corollary}\label{p:hypothesis}
With the notation of Proposition \ref{p::hypothesis}, let $\varepsilon\in (0, 1)$ and $T\in\Bc(\Hc)$ symmetric. Assume that
$$\langle\Lambda\rangle^{\varepsilon}T\in\Bc(\ell^2(\Vc,1)),$$  
then $T\in \Cc^{0,1}(A).$ 
\end{corollary}
\subsubsection{Unitary transformation}
By perturbing the metric, a second Hilbert space appears $\ell^2(\Vc,m)$, which is equal to $\ell^2(\Vc,1)$ but is endowed with a different and equivalent norm. The problem is that $\Delta_m$ is not self-adjoint in $\ell^2(\Vc,1)$. To circumvent this difficulty, we rely on the following transformation: 
\begin{proposition}\label{X}Set the following map  \begin{align*} T_{1\rightarrow m}: \ell^2(\Vc,1)&\to\ell^2(\Vc,m)
 \\
f&\mapsto T_{1\rightarrow m}f(n):=\frac{1}{\sqrt{m(n)}}f(n).\end{align*} Then, the transformation $T_{1\rightarrow m}$ is unitary.
\end{proposition}
\proof Let $f\in\ell^2(\Vc,1)$,
 \begin{align*}  \|T_{1\rightarrow m}f\|_{\ell^2(\Vc,m)}^2&=\sum_{(n_1,n_2)\in\Vc}m(n_1,n_2) \left|T_{1\rightarrow m}f(n_1,n_2)\right|^2\\
&=\sum_{(n_1,n_2)\in\Vc} m(n_1,n_2)\left|\frac{1}{\sqrt{m(n_1,n_2)}}f(n_1,n_2)\right|^2
\\
&=\sum_{(n_1,n_2)\in\Vc} \left|f(n_1,n_2)\right|^2= \|f\|_{\ell^2(\Vc,1)}^2.
\end{align*}
This ensures the result.
\qed

Recalling \eqref{dm} and the hypotheses $(H_0)$. 
Thanks to the unitary transformation, we can transport $\Delta_{m}$ into $\ell^2(\Vc,1)$. Namely,
let $\widetilde{\Delta}:= T^{-1}_{1\rightarrow m}\Delta_{m}T_{1\rightarrow m}$.
\begin{proposition}We have:
 \begin{align}\label{e:transform}
\nonumber\widetilde{\Delta}
=&\frac{1}{6\sqrt{m(Q)}}
 \Big(\frac{1}{\sqrt{m(Q_1+1,Q_2)}}U_1^*\otimes1+\frac{1}{\sqrt{m(Q_1-1,Q_2)}}U_1\otimes1\\
 &+\frac{1}{\sqrt{m(Q_1,Q_2+1)}}1\otimes U_2^*
 +\frac{1}{\sqrt{m(Q_1,Q_2-1)}}1\otimes U_2\\
 &\nonumber+\frac{1}{\sqrt{m(Q_1+1,Q_2-1)}}U_1^*\otimes U_2+\frac{1}{\sqrt{m(Q_1-1,Q_2+1)}}U_1\otimes U_2^*
 \Big).\end{align}
\end{proposition}
We derive the next expression for the perturbation:

\noindent Given $l,n\in\Vc$, we denote by $l\sim n$ if $\Ec(n,l)>0.$ 
\begin{proposition}\label{compact}
\item[1.]For all $f\in\ell^2(\Vc,1)$, we have:
\begin{align}\label{-}
\left(\Delta_T-\widetilde{\Delta}\right)f(n)
\nonumber:=&\frac{1}{6}\Bigg(\left(1-\frac{1}{\sqrt{m(n)}}\right)\sum_{l\sim n}f(l)\\
&\nonumber+\frac{1}{\sqrt{m(n)}}\sum_{l\sim n}\left(1-\frac{1}{\sqrt{m(l)}}\right)f(l)\Bigg)\\
&=(1 - R) \Delta_T + R \Delta_T (1 - R),
\end{align}
where $R(Q):=\frac{1}{\sqrt{m(Q)}}.$
\item [2.]If $(H_0)$ hold true, we have $\Delta_T-\widetilde{\Delta}$ is a compact operator in $\ell^2(\Vc,1).$
\end{proposition}
\proof
\item[1.] This is a straightforward calculus.
\item[2.] We have
\begin{align*}
\left(\Delta_T-\widetilde{\Delta}\right)f(n_1,n_2)
&=\left(1-\frac{1}{\sqrt{m(n)}}\right)\Delta_Tf(n)+\underbrace{\frac{1}{6\sqrt{m(n)}}\sum_{l\sim n}\left(1-\frac{1}{\sqrt{m(l)}}\right)f(l)}_{\widetilde{K_T}f(n)}.
\end{align*}
By using the hypothesis $(H_0)$, we have $1-\frac{1}{\sqrt{m(n)}}\rightarrow 0$, if $n\rightarrow \infty$ and $\Delta_T$ is bounded then $\left(1-\frac{1}{\sqrt{m(\cdot)}}\right)\Delta_T$ is compact. Now, we will show that $\widetilde{K_T}$ is compact. To show that $\widetilde{K_T}$ is compact, it is enough to use that:
$$(\widetilde{K_T}f)(n)=\left(\left(\frac{1}{\sqrt{m(Q)}}\sum_{(j,k)\in\{*,0,1\}^2,j\neq k}U_1^j\otimes U^{k}_2\left(1-\frac{1}{\sqrt{m(Q)}}\right)\right)f\right)(n).$$
In view of the boundedness of $U_1^j\otimes U^{k}_2$ and since $\left(1-\frac{1}{\sqrt{m(\cdot)}}\right)$ is compact 
we have that the operator $\widetilde{K_T}$ is a finite sum of compact operators. So we obtain the result.
\qed
\begin{proposition}\label{com}
Let $m$ and $V$ be two real-valued bounded functions satisfying respectively $(H_0)$ and $(H'_0)$.
We have:
\item[1.]$\Delta_{m}+V(Q)$ is self-adjoint and bounded.
\item[2.]$\sigma_{\rm ess}(\Delta_{m}+V(Q))=\sigma_{\rm ess}(\Delta_{T})$.
\end{proposition}
\proof
 \item[1.] Hypothesis $(H'_0)$ assures the compactness of $V(Q)$. Since $\Delta_{T}$ is self-adjoint and according Theorem \cite[Theorem XIII.14]{RS}, we have that $\Delta_{m}+V(Q)$ is self-adjoint. 
\item[2.] Using the fact that $V(Q)$ is compact and thanks to Proposition \ref{compact}, we deduce that $\sigma_{\rm ess}(\widetilde{\Delta}+V(Q))=\sigma_{\rm ess}(\Delta_{T})$. By using the unitary transformation, we obtain $\sigma_{\rm ess}(\Delta_m+V(Q))=\sigma_{\rm ess}(\Delta_{T})$.\qed
\subsubsection{Perturbed potential}
We start by treating the regularity properties of the potential $V$. The perturbation of the metric will be more involved and treated in the next subsection.
\begin{lemma}\label{H'}
Let $V:\Vc\rightarrow \mathbb{R}$ be a function. We assume that $(H'_0),\ (H'_1),\ (H'_2)$ and $(H'_3)$ hold true, then $V(Q)\in \Cc^{1}(A)$ and $[V(Q),\rmi A]_{\circ}\in \Cc^{0,1}(A)$. In particular, we obtain $V(Q)\in \Cc^{1,1}(A).$
\end{lemma}
\proof
Recalling \eqref{A}. We show this lemma in two steps. First, we prove that $V(Q)\in \Cc^{1}(A)$. It suffices to show that there exists $c>0$, such that: 
$$\left\|\left[V(Q),\rmi A\right]f\right\|^2\leq c\|f\|^2,\ \forall f\in \Sc.$$
Second, we prove that $[V(Q),\rmi A]_{\circ}\in\Cc^{0,1}(A)$. Given $\varepsilon'> 0$, it is enough to show there exists $c_{\varepsilon'}>0$ such that: 
$$\left\|\Lambda^{\varepsilon'}(Q)\left[V(Q),\rmi A\right]f\right\|^2\leq c\|f\|^2,\ \forall f\in \Sc.$$      
Take $\varepsilon'\in [0,1)$ such that $\varepsilon'<\varepsilon$. We work on $\Sc$. We have:
\begin{align*}
\left[V(Q),\rmi A\right]
=&\frac{1}{6}\Big(-\left(Q_1+\frac{1}{2}\right)\left[V(Q),U_1^*\otimes1\right]+\left(Q_1-\frac{1}{2}\right)\left[V(Q),U_1\otimes1\right]\\
&-\left(Q_2+\frac{1}{2}\right)\left[V(Q),1\otimes U_2^*\right]+\left(Q_2-\frac{1}{2}\right)\left[V(Q),1\otimes U_2\right]\\
&-\left(Q_2-Q_1+1\right)\left[V(Q),U_1\otimes U_2^*\right]-\left(Q_1-Q_2+1\right)\left[V(Q),U_1^*\otimes U_2\right]\Big).
\end{align*}
We assume that $(H'_1)$, $(H'_2)$ and $(H'_3)$ are true. Let $f\in\Sc$, we have:
\begin{align}\label{part0}
\nonumber\left\|\Lambda^{\varepsilon'}(Q_1,Q_2)\left[V(Q),\rmi A\right]f\right\|
\nonumber\leq&\frac{1}{6}\Bigg(\left\|\Lambda^{\varepsilon'}(Q_1,Q_2)\left(Q_1+\frac{1}{2}\right)\left[V(Q),U_1^*\otimes1\right]f\right\|\\
\nonumber&+\left\|\Lambda^{\varepsilon'}(Q_1,Q_2)\left(Q_1-\frac{1}{2}\right)\left[V(Q),U_1\otimes1\right]f\right\|\\
\nonumber&+\left\|\Lambda^{\varepsilon'}(Q_1,Q_2)\left(Q_2+\frac{1}{2}\right)\left[V(Q),1\otimes U_2^*\right]f\right\|\\
&+\left\|\Lambda^{\varepsilon'}(Q_1,Q_2)\left(Q_2-\frac{1}{2}\right)\left[V(Q),1\otimes U_2\right]f\right\|\\
\nonumber&+\left\|\Lambda^{\varepsilon'}(Q_1,Q_2)\left(Q_2-Q_1+1\right)\left[V(Q),U_1\otimes U_2^*\right]f\right\|\\
\nonumber&+\left\|\Lambda^{\varepsilon'}(Q_1,Q_2)\left(Q_1-Q_2+1\right)\left[V(Q),U_1^*\otimes U_2\right]f\right\|\Bigg)\\
\nonumber\leq& c_{\varepsilon'} \|f\|.
\end{align}
Here, we have used $(H'_1)$ for the first and second term, $(H'_2)$ for the third and fourth term and $(H'_3)$ for the fifth and sixth term.
Then, taking $\varepsilon' = 0$ by density and thanks to \cite[Lemma 6.2.9]{ABG}, we obtain $V(Q)\in \Cc^{1}(A)$. Next, given $\varepsilon'>0$, with the help of Corollary \ref{p:hypothesis}, \eqref{part0} ensures that $[V(Q),\rmi A]_{\circ}\in \Cc^{0,1}(A)$ and therefore that $V(Q)\in\Cc^{1,1}(A).$ 
\qed
\subsubsection{Perturbed metric}
We turn to the most technical part, the perturbation of the metric and start with a lemma.
\begin{lemma}\label{QQ}
We assume that $(H_0)$, $(H_1)$, $(H_2)$ and  $(H_3)$ hold true, we have $R(Q)\in \Cc^{1}(A)$ and for $\varepsilon' \in [0, \epsilon]$, $\Lambda^{\varepsilon'}(Q) [A,R(Q)]_{\circ}$ is bounded.
\end{lemma}
\proof
To prove that $R(Q)\in \Cc^{1}(A)$. It suffices to show that there exists $c>0$, such that
\begin{equation}\label{eq1}
\left\|\left[R(Q), \rmi A\right]f\right\|^2\leq c \|f\|^2,\ \forall f\in \Sc.
\end{equation}
Let $\sigma \in \{-1, 1\}$.
To simplify, we write $U^{\sigma}_{1}$ for $U^{\sigma}_{1} \otimes 1$ and $U^{\sigma}_{2}$ for $1 \otimes U^{\sigma}_{2}$.
As operators on $\mathcal{S}$, we have, for $j\in \{1,2\}$
\begin{align*}
&\left[(Q_j - \sigma/2) U^{\sigma}_{j}, R(Q_1,Q_2)\right]_{\circ} = \big(Q_j - \sigma/2\big) \left[U^{\sigma}_{j}, R(Q_1,Q_2)\right]_{\circ} \\
&\quad= \big(Q_j - \sigma/2\big) \Big(R(Q_1-\sigma\delta_{j,1},Q_2-\sigma\delta_{j,2}) - R(Q_1,Q_2)\Big) U^{\sigma}_{j},
\end{align*}
where $\delta_{j,i}$ is the Kronecker's delta symbol
and
\begin{align*}
&\left[(\sigma(Q_2 - Q_1) + 1) U^{\sigma}_{1} \otimes U^{-\sigma}_{2}, R(Q_1,Q_2)\right]_{\circ}\\
&= \big(\sigma(Q_2 - Q_1) + 1\big) \left[U_{1}^{\sigma} \otimes U^{-\sigma}_{2}, R(Q_1,Q_2)\right]_{\circ}\\
&= \big(\sigma(Q_2 - Q_1) + 1\big) \Big(R(Q_1-\sigma,Q_2+\sigma) - R(Q_1,Q_2)\Big) U^{\sigma}_{1} \otimes U^{\sigma}_{2}.
\end{align*}
Now,\\
$
|\left(U^{\sigma}_{1}R\right)(n_1, n_2) - R(n_1, n_2)|= \frac{|\eta(n_1 - \sigma, n_2) - \eta(n_1, n_2)|}{\sqrt{m(n_1 - \sigma, n_2)} \sqrt{m(n_1, n_2)} \left(\sqrt{m(n_1 - \sigma, n_2)} + \sqrt{m(n_1, n_2)}\right)},
$\\
thus, from $(H_1)$, we derive that 
\[
M_{1,\sigma} := \sup_{(n_1, n_2) \in \mathbb{Z}^2} \Lambda^\varepsilon(n_1, n_2) \langle n_1 \rangle |\eta(n_1 - \sigma, n_2 ) - \eta(n_1, n_2)| < \infty,
\]
 and since $m \geq c$, for some constant $c > 0$, we infer that
\[
|R(n_1-\sigma, n_2) - R(n_1, n_2)| \leq M_{1,\sigma} \langle n_1 \rangle^{-1} \Lambda^{-\varepsilon}(n_1, n_2) (2c\sqrt{c})^{-1}.
\]
Similarly, from $(H_2)$, we deduce that 
\[
M_{2,\sigma} := \sup_{(n_1, n_2) \in \mathbb{Z}^2} \Lambda^\varepsilon(n_1, n_2) \langle n_2 \rangle |\eta(n_1 , n_2-\sigma ) - \eta(n_1, n_2)| < \infty,
\]
and, as above, we obtain
\[
|R(n_1, n_2-\sigma) - R(n_1, n_2)| \leq M_{2,\sigma} \langle n_2 \rangle^{-1} \Lambda^{-\varepsilon}(n_1, n_2) (2c\sqrt{c})^{-1}.
\]
From $(H_3)$, we conclude that
\[
M_{3,\sigma} := \sup_{(n_1, n_2) \in \mathbb{Z}^2} \Lambda^\varepsilon(n_1, n_2) \langle n_1 - n_2 \rangle |\eta(n_1 - \sigma, n_2 + \sigma) - \eta(n_1, n_2)| < \infty,
\]
and, as above, we get
\[
|R(n_1-\sigma, n_2+\sigma) - R(n_1, n_2)| \leq M_{3,\sigma} \langle n_1 - n_2 \rangle^{-1} \Lambda^{-\varepsilon}(n_1, n_2) (2c\sqrt{c})^{-1}.
\]
Coming back to the commutators, this yields, the boundedness of
\[
\Lambda^{\varepsilon'}(Q_1,Q_2) \left[(Q_j - \sigma/2) U^{\sigma}_{j}, R(Q)\right]_{\circ}\] and \[ \Lambda^{\varepsilon'}(Q) \left[(\sigma(Q_2 - Q_1) + 1) U^{\sigma}_{1} \otimes U^{-\sigma}_{2}, R(Q)\right]_{\circ},
\]
for any $\varepsilon' \in [0, \varepsilon]$. In view of \eqref{A}, this shows that $R(Q) \in \mathcal{C}^1(A)$ and $\Lambda^{\varepsilon'} [A,R(Q)]_{\circ}$ is bounded.\qed
\begin{proposition}\label{H}
We assume that $(H_0)$, $(H_1)$, $(H_2)$ and  $(H_3)$ hold true, we have $\widetilde{\Delta}\in \Cc^{1}(A)$. Moreover $[\widetilde{\Delta},\rmi A]_{\circ}\in\Cc^{0,1}(A)$. In particular, $\widetilde{\Delta}\in \Cc^{1,1}(A).$
\end{proposition}
\proof 
First of all, we recall \eqref{A} and \eqref{-}. The proof is constituted as follows:
In the first step, we are going to prove that $\widetilde{\Delta}\in \Cc^{1}(A)$. It suffices to show that there exists $c>0$, such that:
\begin{equation}\label{eq1}
\left\|\left[\Delta_T-\widetilde{\Delta}, \rmi A\right]f\right\|^2\leq c \|f\|^2,\ \forall f\in \Sc.
\end{equation}
Then, by density and thanks to \cite[Proposition 6.2.9]{ABG} and Proposition \ref{C2}, we obtain the result.
In the second step, we will establish that $[\widetilde{\Delta},\rmi A]_{\circ}\in\Cc^{0,1}(A)$. Given $\varepsilon'\in[0,\varepsilon],\ \varepsilon\in (0,1)$, we show there
exists $c_{\varepsilon'} > 0$ such that:
\begin{equation}\label{eq2}
\left\|\Lambda^{\varepsilon'}(Q_1,Q_2)\left[\Delta_T-\widetilde{\Delta},\rmi A\right]f\right\|^2\leq c_{\varepsilon'}\|f\|^2,\ \forall f\in \Sc.
\end{equation}
Then, by density, $[\Delta_T-\widetilde{\Delta},\rmi A]_{\circ}\in \Cc^{0,1}(A)
$. Finally, by
Corollary \ref{p:hypothesis} and thanks to Proposition \ref{C2}, we have $[\widetilde{\Delta},\rmi A]_{\circ}\in \Cc^{0,1}(A)$. In particular,  thanks to Proposition \ref{C2}, we obtain $\widetilde{\Delta}\in \Cc^{1,1}(A).$\\
Thus, as operators acting on $\mathcal{S}$, due to simplifications, we use \eqref{-}, we obtain 
\begin{align}\label{1}
&[A, \Delta_T - \tilde{\Delta}] \\
&\nonumber= [A, \Delta_T] - R(Q)[A, \Delta_T]R(Q) - [A, R(Q)]\Delta_T R(Q) - R(Q)\Delta_T [A, R(Q)]. 
\end{align}
By Lemma \ref{QQ}, by Proposition \ref{C2}, we obtain that $R(Q)\in\Cc^1(A)$, then we know that the closure of first two terms on the r.h.s of \eqref{1} are in $\mathcal{C}^1(A)$, since $\Delta_T \in \mathcal{C}^2(A)$.
Moreover, for $\varepsilon' \in [0, \varepsilon]$, $\Lambda^{\varepsilon'}(Q) [A, R(Q)]_{\circ} \Delta_T R(Q)$ is bounded, since $R(Q)$ and $\Delta_T$ are bounded and by Lemma \ref{QQ}. Then $[A, R(Q)]_{\circ}\Delta_T R(Q)\in\Cc^{0,1}(A).$\\
Next, we prove that $\Lambda^{\varepsilon'}(Q) \Delta_T \Lambda^{-\varepsilon'}(Q)$ is bounded.
We obtain the boundedness of $\Lambda(Q) \Delta_T \Lambda^{-1}(Q)$ from 
\[
\Lambda(Q) \Delta_T \Lambda^{-1}(Q) = \Delta_T + [\Lambda(Q), \Delta_T]_{\circ} \Lambda^{-1}(Q)
\]
and a direct computation of the commutator of $\Lambda$ with the $U^{\sigma}_{1} \otimes 1$, $1 \otimes U^{\sigma}_{2}$, and $U^{\sigma}_{1} \otimes U^{-\sigma}_{2}$,
at the end, we conclude by interpolation, as in Remark \ref{Lambdavarepsilon}.\\
Finally, thanks to  Lemma \ref{QQ}, we obtain $\Lambda^{\varepsilon'}(Q) [A, R(Q)]_{\circ}$ is bounded and we have $\Lambda^{\varepsilon'}(Q) \Delta_T \Lambda^{-\varepsilon'}(Q)$ is bounded as well, then we see that $$\Lambda^{\varepsilon'}(Q) R(Q) \Delta_T [A, R(Q)]_{\circ}= R(Q)\Lambda^{\varepsilon'}(Q) \Delta_T\Lambda^{-\varepsilon'}(Q)\Lambda^{\varepsilon'}(Q) [A, R(Q)]_{\circ}$$ is bounded. Thus $\Delta_T - \tilde{\Delta}$ is $\mathcal{C}^{1}(A)$ by taking $\varepsilon'=0$ and $\mathcal{C}^{1,1}(A)$ considering $\varepsilon'\in ]0,\varepsilon].$
\qed
\subsection{Proof of the main result}\label{P.result}
The main result of this section is Theorem \ref{t:LAP}.
To begin with, we establish the Mourre estimate in the case of perturbation. As the ambient space is now $\ell^2(\Vc,m)$, we transport the operators acting in $\ell^2(\Vc,m)$ into it. We start with a remark.
\begin{remark}\label{XXXX}
Recalling Proposition \ref{X}. Let $A_m:=T_{1\rightarrow m}AT^{-1}_{1\rightarrow m}, \forall z\in\mathbb{C}\setminus\mathbb{R},$ $$T_{1\rightarrow m}(A-z)^{-1}T^{-1}_{1\rightarrow m}=(A_m-z)^{-1}.$$ By functional calculus this gives $T_{1\rightarrow m}e^{\rmi tA}T^{-1}_{1\rightarrow m}=e^{\rmi tA_m}$. In turn given $S$ bounded in $\ell^2(\Vc,1)$, we have $S\in \Cc^{\alpha}(A)\Leftrightarrow T_{1\rightarrow m}ST^{-1}_{1\rightarrow m}\in \Cc^{\alpha}(A_m)$, with $\alpha\in \{$1$;$2$;$0,1$;$1,1$\}$ and for $\alpha=1,\ [T_{1\rightarrow m}ST^{-1}_{1\rightarrow m},\rmi A_m]_{\circ}=T_{1\rightarrow m}[S,\rmi A]_{\circ}T^{-1}_{1\rightarrow m}.$
\end{remark}
\noindent Next, since $V(Q)$ is an operator of multiplication, so $V(Q):=T_{1\rightarrow m}V(Q)T^{-1}_{1\rightarrow m}$. Consequently, we have $$T_{1\rightarrow m}(\widetilde{\Delta}+V(Q))T_{1\rightarrow m}^{-1}=T_{1\rightarrow m}\widetilde{\Delta}T^{-1}_{1\rightarrow m}+T_{1\rightarrow m}V(Q)T^{-1}_{1\rightarrow m}=\Delta_m+V(Q):=H_m.$$  
\begin{theorem}\label{mou}
Let $V:\Vc\rightarrow \mathbb{R}$. We assume that $(H_0),\ (H_1),\ (H_2)$, $(H_3),\ (H'_0),$ $(H'_1),$ $(H'_2)$ and $(H'_3) $ hold true.
Then $\widetilde{\Delta}+V(Q)\in \Cc^{1,1}(A)$. Moreover, for all compact interval $\Ic\subset[-\frac{1}{2},1]\backslash\{1,-\frac{1}{2},-\frac{1}{3}\}$, there are $c>0$ and a compact operator $\widetilde{K}$ such that:
\begin{align}\nonumber
&E_\Ic(\widetilde{\Delta}+V(Q))\left[\widetilde{\Delta}+V(Q), \rmi A\right]_\circ E_\Ic(\widetilde{\Delta}+V(Q))\\
&\label{Mour}\hspace*{+3.8cm}\geq c E_\Ic(\widetilde{\Delta}+V(Q))+K.
\end{align}
Equivalently, $H_m\in\Cc^{1,1}(A_m)$ and  
\begin{align}
&E_\Ic(H_m)\left[H_m, \rmi A_m\right]_\circ E_\Ic(H_m)\label{Mou}\geq c E_\Ic(H_m)+K_m,
\end{align}
where $K_m:=T_{1\rightarrow m}KT^{-1}_{1\rightarrow m}.$
\end{theorem}
\proof 
The Proposition \ref{H}, the Lemma \ref{H'} and the Proposition \ref{C2} give that $\widetilde{\Delta}+V(Q)\in \Cc^{1,1}(A)$. By hypotheses $V(Q)$ is compact and by Proposition \ref{compact}, we have $\left(\Delta_T-\widetilde{\Delta}\right)$ is a compact operator. Then, by using Proposition \ref{Mourre} and by \cite[Theorem 7.2.9]{ABG}, we obtain \eqref{Mour}. Using the transformation unitary $T_{1\rightarrow m}$, Remark \ref{XXXX} ensures \eqref{Mou}.\qed\\
\textbf{Proof of Theorem \ref{t:LAP}}:
Proposition \ref{com} provides point $1.$ and Theorem \ref{mou} gives the points $2.$  
To show point $4.$ it is enough to consider $s>\frac{1}{2}$. We apply \cite[Proposition 7.5.6]{ABG} and we obtain:
$$\lim_{\rho\rightarrow 0^{+}}\sup_{\lambda\in[a,b]}\|\langle A_m\rangle^{-s}(H_m-\lambda-\rmi\rho)^{-1}\langle A_m\rangle^{-s}\| \mbox{ is finite.}$$
Furthermore, in the norm topology of bounded operators, the boundary values of the
resolvent:
\[ [a,b] \ni\lambda\mapsto\lim_{\rho\to0^{\pm}}\langle A_m\rangle^{-s}(H_m-\lambda-\rmi\rho)^{-1}\langle A_m\rangle^{-s} \mbox{ exists and is continuous},\]
where $[a,b]$ is included in $\mathbb{R}\setminus\left(\kappa(H_m)\cup\sigma_{\rm p}(H)\right)$. In particular, this gives Point $3.$
 By Remark \ref{Lambdavarepsilon}, there is $c>0$ such that: $$\|\langle A_m\rangle^{s}f\|\leq c\|\Lambda^{s}(Q)f\|,$$
for all $f\in\Dc(\Lambda^{s}(Q)).$ We conclude that $$\lim_{\rho\rightarrow 0^{+}}\sup_{\lambda\in[a,b]}\| \Lambda^{-s}(Q)(H_m-\lambda-\rmi\rho)^{-1}\Lambda^{-s}(Q)\| \mbox{ is finite.}$$
 The point $5.$ is a consequence of $4.$\qed

\end{document}